\documentclass{amsart}

\usepackage[colorlinks=true, urlcolor=black, citecolor=black, linkcolor=black, hyperfootnotes=true]{hyperref}
\usepackage{amssymb}
\usepackage{aliascnt}
\usepackage{mathscinet}
\usepackage{graphicx}
\usepackage{accents}
\usepackage{enumitem}

\numberwithin{equation}{section}

\newtheorem{thm}{Theorem}[section]

\newaliascnt{prp}{thm}
\newtheorem{prp}[prp]{Proposition}
\aliascntresetthe{prp}

\newaliascnt{cor}{thm}
\newtheorem{cor}[cor]{Corollary}
\aliascntresetthe{cor}

\theoremstyle{definition}

\newaliascnt{dfn}{thm}
\newtheorem{dfn}[dfn]{Definition}
\aliascntresetthe{dfn}

\newtheorem*{qst}{Question}

\newaliascnt{rmk}{thm}
\newtheorem{rmk}[rmk]{Remark}
\aliascntresetthe{rmk}

\author{Tristan Bice}
\email{tristan.bice@gmail.com}
\thanks{The author is supported by the GA\v{C}R project EXPRO 20-31529X and RVO: 67985840 at the Institute of Mathematics of the Czech Academy of Sciences in Prague, Czech Republic}

\title{Gr\"atzer-Hofmann-Lawson-Jung-S\"underhauf Duality}

\begin{document}

\begin{abstract}
We unify several extensions of the classic Stone duality due to Gr\"atzer, Hoffman-Lawson and Jung-S\"underhauf.  Specifically we show that $\cup$-bases of locally compact sober spaces are dual to $\prec$-distributive $\vee$-predomains, where $\prec$ is a transitive relation representing compact containment.
\end{abstract}

\maketitle

\section*{Introduction}

\subsection*{Historical Background}

Over 80 years ago in \cite{Stone1936}, Stone exhibited a remarkable duality between Boolean algebras and clopen bases of compact Hausdorff spaces, one  which inspired a whole host of dualities between algebraic and more analytic structures.  Even the specific duality that Stone examined would be extended several times in the ensuing decades.  The first step was taken by Stone himself in \cite{Stone1938}, where he extended to the duality from Boolean algebras to general distributive lattices, and from clopen bases of compact Hausdorff spaces to compact open bases of locally compact sober spaces that are coherent.  It was only much later in \cite[II.5]{Gratzer1978} that coherence was eliminated by Gr\"atzer, the key idea being to work with $\vee$-semilattices satisfying an appropriate generalisation of distributivity.

Around the same time, Hofmann and Lawson obtained a duality for general locally compact sober spaces in \cite{HofmannLawson1978}, but at the cost of dealing with big lattices, namely continuous frames representing the entire open set lattice.  It was only much later again in \cite{JungSunderhauf1995}, building on ideas from \cite{Smyth1992}, that a partial unification of the Hofmann-Lawson and Stone dualities was achieved.  However, there were still a couple of issues, e.g. Jung-S\"underhauf's work was again restricted to coherent spaces.  Also, their axioms were too strong to apply to the entire open set lattice, and even finding an appropriate basis to satisfy their axioms turned out to be a non-trivial task.

Nevertheless, the basic idea of \cite{JungSunderhauf1995} was perfectly sound, namely to work with a lattice, or better yet a $\vee$-semilattice, together with an extra transitive relation $\prec$ representing compact containment.  However, instead of imposing proximity-like conditions on $\prec$, we take leaf out of Gr\"atzer's book and formulate an even more general notion of distributivity for $\prec$.  The theory then proceeds more smoothly along classical lines (e.g. extending the Birkhoff-Stone prime filter theorem to $\prec$), allowing us to simultaneously generalise and hence unify all these Stone duality extensions.  Moreover, we make the duality functorial via certain relational morphisms, just like in \cite{JungSunderhauf1995} and \cite{BiceStarling2018}.

\subsection*{Related Work}

In doing this, we bring Stone's duality much closer to a contemporaneous duality by Wallman in \cite{Wallman1938}.  Indeed, it is an interesting artefact of history that Wallman's duality was largely overlooked in favour of Stone's duality, even though Wallman's duality applies to (e.g. connected) spaces more commonly seen in other fields.  However, even the most general form of Stone duality presented here differs in some important respects to Wallman's duality.  For example, the former deals with sober spaces while the latter deals with $T_1$ spaces, which are incomparable notions even among compact spaces (see \cite[I.3.1]{PicadoPultr2012}).  Also Wallman's duality applies to compact spaces, even those that are not locally compact (in the sense of points having compact neighbourhood bases).  Recently in \cite{BiceKubis2020}, we extended Wallman's duality to `locally closed compact' spaces which, in the Hausdorff case, coincide with locally compact spaces.  But even in the this case the dualities differ in that Wallman's duality can be further generalised to semilattices representing mere subbases (again see \cite{BiceKubis2020}).

In a different direction, the original Stone duality in \cite{Stone1936} can be generalised by reducing the reliance on joins/unions.  This is particularly important for obtaining non-commutative extensions to inverse semigroups and \'etale groupoids.  For example, in \cite{BiceStarling2018} we showed how to extend the classic Stone duality in \cite{Stone1936} to locally Hausdorff spaces (even non-zero dimensional ones) by considering conditional $\vee$-semilattices.  Even without any joins whatsoever, a duality can be still obtained with bases or even pseudobases of locally compact Hausdorff spaces.  This we showed in \cite{BiceStarling2019}, providing a common generalisation of Exel's tight groupoid construction from \cite{Exel2008} and the earlier Stone duality extensions of Shirota \cite{Shirota1952} and de Vries \cite{deVries1962}.  Note, however, that without joins we are forced to rely more heavily on stronger separation properties, like being (locally) Hausdorff.

\subsection*{Outline}

First in \autoref{TheSpectrum}, we examine the spectrum of proper round prime filters in any $\vee$-semilattice $S$ with minimum $0$ together with an additional relation $\prec$.  The key results are that the spectrum is always sober \textendash\, see \autoref{RoundSpectrumSober} \textendash\, and that the spectrum recovers any core compact sober space from any $\cup$-basis \textendash\, see \autoref{SoberRecovery}.

Next, in \autoref{DistributivePredomains}, we investigate various properties of $\prec$ and their relation to the spectrum.  Most of these properties have been considered before in domain theory \textendash\, see \cite{GierzHofmannKeimelLawsonMisloveScott2003}.  The key extra condition we consider is a generalisation of distributivity for $\prec$, a variant of a previous version in \cite{BiceStarling2018} \textendash\, see \eqref{precDistributive}.  This allows for a direct generalisation of the Birkhoff-Stone prime filter theorem in \autoref{BirkhoffStone}.  Together with the conditions defining $\vee$-predomains, it then follows that the spectrum is core compact and that $\leq$ and $\prec$ are faithfully represented on the spectrum as $\subseteq$ and $\Subset$ respectively \textendash\, see \autoref{CoreCompactSpectrum}, \eqref{leq=>subseteq'} and \eqref{prec=>Subset'}.

Finally, in \autoref{MorphismsS}, we make the duality functorial with respect to certain relational morphisms defined by a few simple first order properties \textendash\, see \autoref{Morphisms}.  With a couple of extra conditions on the relational morphisms, this even yields an equivalence between categories of $\prec$-distributive $\vee$-predomains and $\cup$-bases of locally compact sober spaces.

\newpage

\section{The Spectrum}\label{TheSpectrum}

We make the following standing assumption throughout the paper.
\[\textbf{We have a relation $\prec$ on a $\vee$-semilattice $S$ with minimum }0.\]

First we recall the following standard definitions from order theory.
\begin{dfn}
We call $F\subseteq S$ a \emph{filter} and $I\subseteq S$ an \emph{ideal} if
\begin{align}
\tag{Filter}p,q\in F\qquad&\Leftrightarrow\qquad\exists r\in F\ (r\leq p,q).\\
\tag{Ideal}p,q\in I\qquad&\Leftrightarrow\qquad\exists r\in I\ (r\geq p,q).\\
\intertext{A filter $P\subseteq S$ is \emph{prime} if $S\setminus P$ is an ideal.  As $S$ is a $\vee$-semilattice, this means}
\tag{Prime}p\vee q\in P\qquad&\Rightarrow\qquad p\in P\text{ or }q\in P.\\
\intertext{A filter $R\subseteq S$ is \emph{round} if}
\tag{Round}p\in R\qquad&\Rightarrow\qquad\exists r\in R\ (r\prec p).
\end{align}
\end{dfn}

\begin{dfn}
The \emph{spectrum} of $(S,\leq,\prec)$ is the space with points
\[\widehat{S}=\{P\subseteq S:P\text{ is a proper round prime filter}\}\]
with the topology generated by $(\widehat{S}_p)_{p\in S}$ where
\[\widehat{S}_p=\{P\in\widehat{S}:p\in P\}.\]
\end{dfn}

Note that, as $\widehat{S}$ consists of filters, $(\widehat{S}_p)_{p\in S}$ is a basis for the topology on $\widehat{S}$.

\begin{dfn}
A closed set $C\subseteq X$ is \emph{irreducible} if its proper closed subsets form an ideal.  A space $X$ is \emph{sober} if every irreducible $C$ has a unique dense point $x$, i.e.
\[C=\mathrm{cl}\{x\}.\]
\end{dfn}

\begin{prp}\label{RoundSpectrumSober}
The spectrum is always sober.
\end{prp}

\begin{proof}
Take any irreducible $C\subseteq\widehat{S}$ and let
\[P=\{p\in S:\widehat{S}_p\cap C\neq\emptyset\}.\]
If $q\geq p\in P$ then $\widehat{S}_q\cap C\supseteq\widehat{S}_p\cap C\neq\emptyset$ and hence $q\in P$ too.  Also, for any $p,q\in P$, $\widehat{S}_p\cap C\neq\emptyset\neq\widehat{S}_q\cap C$ implies $\widehat{S}_p\cap\widehat{S}_q\cap C\neq\emptyset$, as $C$ is irreducible.  Taking any $Q\in\widehat{S}_p\cap\widehat{S}_q\cap C$, we see that $p,q\in Q$ and hence $r\in Q$, for some $r\leq p,q$, as $Q$ is a filter.  Thus $Q\in\widehat{S}_r\cap C\neq\emptyset$ and hence $r\in P$, showing that $P$ is a filter.  Similarly, for any $p\in P$, we have $Q\in\widehat{S}_p\cap C$ so $p\in Q$ and hence we have $q\in Q$ with $q\prec p$, as $Q$ is round.  Thus $Q\in\widehat{O}_q\cap C\neq\emptyset$ so $q\in P$, showing that $P$ is also round.  For any $p,q\in S\setminus P$, we have $\widehat{S}_{p\vee q}=\widehat{S}_p\cup\widehat{S}_q\subseteq X\setminus C$ so $S\setminus P$ is an ideal.  Also $\widehat{S}_0=\emptyset$ so $0\notin P$ and hence $P$ is a proper round prime filter, i.e. $P\in\widehat{S}$.

As $P\notin\widehat{S}_q$ implies $q\notin P$ and hence $\widehat{S}_q\cap C=\emptyset$, it follows that $C\subseteq\mathrm{cl}\{P\}$.  On the other hand, as $\widehat{S}_p\cap C\neq\emptyset$, for all $p\in P$, we have $P\in\mathrm{cl}(C)=C$ and hence $P$ is a dense point in $C$.  Moreover, $P$ is unique, as $\widehat{S}$ is immediately seen to be $T_0$.
\end{proof}

Conversely, any `core compact' sober space arises in this way.  Specifically, for any topology $\mathcal{O}(X)$ on $X$, let $\Subset$ denote the way-below relation on $\mathcal{O}(X)$, i.e. $p\Subset q$ means that every open cover of $q$ has a finite subcover of $p$:
\[p\Subset q\qquad\Leftrightarrow\qquad\forall\mathcal{C}\subseteq\mathcal{O}(X)\ (q\subseteq\bigcup\mathcal{C}\ \Rightarrow\ \exists\text{ finite }\mathcal{F}\subseteq\mathcal{C}\ (p\subseteq\bigcup\mathcal{F})).\]
As in \cite[\S5.2.1]{Goubault2013}, we call $X$ \emph{core compact} if $\mathcal{O}(X)$ is a continuous frame, i.e. if every open neighbourhood filter $\mathcal{O}_x=\{O\in\mathcal{O}(X):x\in O\}$ is round w.r.t. $\Subset$.

\begin{dfn}
A \emph{$\cup$-basis} $S\subseteq\mathcal{O}(X)$ is a basis closed under finite unions.
\end{dfn}

We include the empty union here, i.e. the empty set $\emptyset=\bigcup\emptyset$.  So, equivalently, a $\cup$-basis is a basis containing $\emptyset$ which is closed under pairwise unions.

\begin{thm}\label{SoberRecovery}
For any $\cup$-basis $S$ of a core compact sober space $X$,
\[x\mapsto S_x=\{s\in S:x\in s\}\]
is a homemorphism from $X$ onto $\widehat{S}$, where $(S,\leq,\prec)=(S,\subseteq,\Subset)$.
\end{thm}

\begin{proof}
Each $S_x$ is certainly a proper prime filter (proper because $\emptyset\in S\setminus S_x$).  As $X$ is core compact, each $S_x$ is also round and hence $S_x\in\widehat{S}$.

Conversely, say $P\in\widehat{S}$.  Let $O=\bigcup(S\setminus P)$ and note that $p\nsubseteq O$, for all $p\in P$ \textendash\, if we had $p\subseteq O$ then, as $P$ is round, we would have $q\in P$ with $q\Subset p\subseteq O$ and hence $q\subseteq\bigcup F$, for some finite $F\subseteq S\setminus P$, contradicting the fact $P$ is prime.  For any open $N\supsetneqq O=\bigcup(S\setminus P)$, we have $p\in P$ with $p\subseteq N$, as $S$ is a basis.  For any other $N'\supsetneqq O$, we again have $p'\in P$ with $p'\subseteq N'$.  As $P$ is a filter, we have $p''\in P$ with $p''\subseteq p\cap p'\subseteq N\cap N'$.  Thus $N\cap N'\supseteq p''\nsubseteq O$, as $p''\in P$, so $O\neq N\cap N'$, showing that $X\setminus O$ is irreducible.  As $X$ is sober, we have $x\in X$ with $X\setminus O=\mathrm{cl}\{x\}$.  As $x\notin O=\bigcup(S\setminus P)$, $S_x\subseteq P$.  On the other hand, if $s\in S\setminus S_x$ then $s\cap\mathrm{cl}\{x\}=\emptyset$ and hence $s\subseteq O$, which can only mean $s\notin P$.  Thus $S_x=P$.

So $x\mapsto S_x$ is a bijection from $X$ to $\widehat{S}$.  Consequently, it is a homeomorphism, as we immediately see that it maps the basis $S$ onto the basis $(\widehat{S}_p)_{p\in S}$.
\end{proof}

Our next goal is to isolate the abstract properties of $(S,\leq,\prec)$ that make $\widehat{S}$ core compact and ensure that $\prec$ and $\leq$ are faithfully represented as $\subseteq$ and $\Subset$ on $(\widehat{S}_p)_{p\in S}$.

\begin{rmk}
A more standard approach would be to take $(S,\leq)$ as the official structure and derive $\prec$ in some way from $\leq$, e.g. if $(S,\leq,\prec)=(\mathcal{O}(X),\subseteq,\Subset)$ then $\prec$ is, by definition, the way-below relation derived from $\leq$ on $S$.  However, we are more interested in $\cup$-bases $S\subseteq\mathcal{O}(X)$.  In this case, as long as $X$ is Hausdorff and each $s\in S$ is relatively compact, $\prec$ can instead be derived as the rather below relation from $\leq$ (see \cite{BiceStarling2018}).  But in more general sober spaces, even the compact and core compact ones, there may actually be no way of defining $\prec$ from $\leq$.

To see this, consider $S=\omega+2=\{0,1,\ldots,\omega,\omega+1\}$ with its usual ordering $\leq$.  On the one hand, we could take $\prec\ =\ \leq$ too and then $(S,\leq,\prec)$ would represent the $\cup$-basis $(\omega+3)\setminus\{\omega\}$ of $X=\omega+2$ in its lower topology.  On the other hand, if we modify $\prec$ just a little by declaring that $\omega\not\prec\omega$ then $(S,\leq,\prec)$ would represent the entire open set lattice of $X=\omega+1$ in its lower topology.

A better idea for our work would be to instead take $(S,\prec)$ as the official structure and derive $\leq$ as the canonical lower preorder of $\prec$ \textendash\, see \eqref{LowerPreorder} below.  We will do this eventually, but initially it suffices to consider two separate relations $\prec$ and $\leq$ that are simply related in various ways.
\end{rmk}

\section{$\prec$-Distributive $\vee$-Predomains}\label{DistributivePredomains}

\begin{dfn}
We call $\prec$ \emph{distributive} if, for all $p,s,t\in S$,
\[\tag{Distributivity}\label{precDistributive}p\prec s\vee t\qquad\Rightarrow\qquad\forall p'\prec p\ \exists s'\prec s\ \exists t'\prec t\ (p'\leq s'\vee t'\leq p).\]
\end{dfn}

Note that when $\prec\ =\ \leq$, it suffices to take $p'=p$ above, which then reduces to the usual notion of distributivity for $\vee$-semilattices.

The motivating example we have in mind is $\Subset$.

\begin{prp}\label{CoreCompactDistributivity}
If $S$ is a $\cup$-basis of a core compact space $X$ then $\Subset$ is distributive.
\end{prp}

\begin{proof}
Say $p'\Subset p\subseteq s\cup t$.  For every $x\in p$, we have $x\in s$ or $x\in t$.  As $S$ is a basis and $X$ is core compact, this means we have $q\in S$ with either $x\in q\Subset p\cap s$ or $x\in q\Subset p\cap t$.  As $p'\Subset p$, finitely many such $q$ cover $p'$.  As $S$ is $\cup$-closed, the union of those $q$ contained in $s$ yields $s'\Subset s$ and the union of those $q$ contained in $t$ yields $t'\Subset t$ such that $p'\subseteq s'\cup t'\subseteq p$.
\end{proof}

Another important property of $\Subset$ is auxiliarity.

\begin{dfn}
We call $\prec$ \emph{auxiliary} if, for all $p,p',q,q'\in S$,
\[\tag{Auxiliarity}\label{Auxiliary}p\leq p'\prec q'\leq q\qquad\Rightarrow\qquad p\prec q\qquad\Rightarrow\qquad p\leq q.\]
\end{dfn}

Note that when $\prec$ is auxiliary, round filters are precisely the \emph{$\prec$-filters}, i.e.
\[\tag{$\prec$-Filter}p,q\in F\qquad\Leftrightarrow\qquad\exists r\in F\ (r\prec p,q),\]
which means that $F$ is $\prec$-directed and upwards closed w.r.t. $\prec$, i.e.
\[\tag{Up-Set}F\subseteq F^\prec=\{p\in S:\exists f\in F\ (f\prec p)\}.\]
When $\prec$ is also distributive, we can extend the Birkhoff-Stone prime filter theorem.

\begin{thm}\label{BirkhoffStone}
If $\prec$ is distributive and auxiliary, $I\subseteq S$ is an ideal and $F\subseteq S$ is a round filter with $I\cap F=\emptyset$ then $F$ extends to a round prime filter $P$ with $I\cap P=\emptyset$.
\end{thm}

\begin{proof}
By Kuratowski-Zorn, we can extend $I$ to an ideal $J$ that is maximal among ideals disjoint from $F$.  Then $P=S\setminus J$ is an up-set and we claim that $P$ is also a round filter.  To see this, take $p,q\in P$.  By maximality, we have $j,k\in J$ with $p\vee j\in F$ and $q\vee k\in F$.  As $F$ is a round filter, we have $f_0,\cdots,f_4\in F$ with $f_4\prec\cdots\prec f_0\leq p\vee j,q\vee k$.  As $\prec$ is auxiliary, $f_1\prec p\vee j,q\vee k$.  As $\prec$ is distributive, we have $p'\prec p$ and $j'\prec j$ with $f_2\leq p'\vee j'\leq f_1$.  Applying auxiliarity and distributivity again, we get $p''\prec p'$ and $j''\prec j'$ with $f_4\leq p''\vee j''\leq f_3$.  As $p''\prec p'\leq f_1\prec q\vee k$, distributivity yet again yields $q'\prec q$ and $k'\prec k$ with $p''\leq q'\vee k'\leq p'$.  Note $q'\prec q$ and $q'\leq p'\prec p$.  Also $j\vee k\in J$ and $j\vee k\vee q'\geq j'\vee k'\vee q'\geq j''\vee p''\geq f_4\in F$ so we must have $q'\in P$.  This shows that $P$ is a $\prec$-filter and hence a round filter.  As $J=S\setminus P$ is an ideal, $P$ is also prime.
\end{proof}

\begin{dfn}
We say $\prec$ is \emph{interpolative} if, for all $p,q\in S$,
\[\tag{Interpolation}p\prec q\qquad\Rightarrow\qquad\exists s\ (p\prec s\prec q).\]
\end{dfn}

Again, the key example to keep in mind is $\Subset$.

\begin{prp}\label{CoreCompactInterpolation}
If $S$ is a $\cup$-basis of core compact $X$ then $\Subset$ has interpolation.
\end{prp}

\begin{proof}
Say $p,q\in S$ satisfy $p\Subset q$.  As $X$ is core compact, for any $x\in q$, we have $r,s\in S$ with $x\in r\Subset s\Subset q$.  As $p\Subset q$, finitely many such $r$ cover $p$.  Taking the union of these $r$'s and $s$'s yields $r',s'\in S$ with $p\subseteq r'\Subset s'\Subset q$, as required.
\end{proof}

\begin{prp}\label{prec=>CompactContainment}
If $\prec$ is distributive, auxiliary and interpolative then
\[p\prec q\qquad\Rightarrow\qquad\exists\text{ compact }C\subseteq\widehat{S}\ (\widehat{S}_p\subseteq C\subseteq\widehat{S}_q).\]
\end{prp}

\begin{proof}
By interpolation, we have $\prec$-filter $F\subseteq p^\prec$ containing $q$ (take the upwards closure of $(p_n)$ where $p\prec p_{k+1}\prec p_k\prec q$, for all $k$).  We claim that $C=\bigcap_{f\in F}\widehat{S}_f$ is compact, which proves the result, as $\widehat{S}_p\subseteq C\subseteq\widehat{S}_q$.  For this it suffices to show
\[\forall j\in I\ (C\nsubseteq\widehat{S}_j)\qquad\Rightarrow\qquad C\nsubseteq\bigcup_{j\in I}\widehat{S}_j,\]
for any non-empty ideal $I\subseteq S$.  But $C\nsubseteq\widehat{S}_j$ implies $j\notin F$, so if this holds for all $j\in I$ then $I\cap F=\emptyset$.  Then \autoref{BirkhoffStone} yields $P\in\widehat{S}$ extending $F$ disjoint from $I$ so $P\in C\setminus\bigcup_{j\in I}\widehat{S}_j$ witnesses $C\nsubseteq\bigcup_{j\in I}\widehat{S}_j$, as required.
\end{proof}

\begin{cor}\label{CoreCompactSpectrum}
If $\prec$ is distributive, auxiliary and interpolative, $\widehat{S}$ is core compact.
\end{cor}

\begin{proof}
If $P\in\widehat{S}_p$ then $q\prec p$, for some $q\in P$.  By \autoref{prec=>CompactContainment}, $P\in\widehat{S}_q\Subset\widehat{S}_p$.
\end{proof}

Actually, what this really shows is that every point in $\widehat{S}$ has a compact neighbourhood base.  By \autoref{SoberRecovery}, together with \autoref{CoreCompactDistributivity} and \autoref{CoreCompactInterpolation}, this provides an alternative proof that core compact sober spaces are necessarily locally compact (see \cite[Theorem 8.3.10]{Goubault2013}).

So \autoref{CoreCompactSpectrum} tells us that the spectrum produces the spaces we want.  The remaining question is what additional conditions ensure that $\subseteq$ and $\Subset$ on the basis $(\widehat{S}_p)_{p\in S}$ defined by $S$ correspond precisely to $\prec$ and $\leq$ on $S$ itself.

\begin{dfn}
We call $\prec$ \emph{approximating} if $\leq$ is the lower preorder of $\prec$, i.e.
\[\tag{Lower Preorder}\label{LowerPreorder}p\leq q\qquad\Leftrightarrow\qquad p^\succ\subseteq q^\succ.\]
\end{dfn}

Note that $\Subset$ is approximating on any basis of a core compact space.

\begin{rmk}
We could have also used the original version of distributivity
\[p\leq s\vee t\qquad\Leftrightarrow\qquad\forall p'\prec p\ \exists s'\prec s\ \exists t'\prec t\ (p'\prec s'\vee t'\prec p)\]
from \cite{BiceStarling2018}, which already implies that $\prec$ is approximating and has interpolation.  We avoided this in order to show that \autoref{BirkhoffStone} only requires the weaker version in \eqref{precDistributive} above.
\end{rmk}

\begin{prp}\label{leq=>subseteq}
If $\prec$ is distributive, auxiliary, interpolative and approximating,
\begin{equation}\label{leq=>subseteq'}
p\leq q\qquad\Leftrightarrow\qquad\widehat{S}_p\subseteq\widehat{S}_q.
\end{equation}
\end{prp}

\begin{proof}
The $\Rightarrow$ part is immediate.  Conversely, say $p\nleq q$.  As $\prec$ is approximating, we have $r\prec p$ with $r\not\prec q$.  As $\prec$ has interpolation, we have a $\prec$-filter $F\subseteq r^\prec$ containing $p$.  As $r\not\prec q$, $F$ is disjoint from the ideal $q^\geq$.  By \autoref{BirkhoffStone}, $F$ extends to a round prime filter $P$ avoiding $q$ so $P\in\widehat{S}_p\setminus\widehat{S}_q$ witnesses $\widehat{S}_p\nsubseteq\widehat{S}_q$.
\end{proof}

\begin{dfn}
We say $\prec$ is $\vee$-preserving if, for all $p,p',q,q'\in S$,
\[p'\prec p\quad\text{and}\quad q'\prec q\qquad\Rightarrow\qquad p'\vee q'\prec p\vee q.\]
$S$ is a \emph{$\vee$-predomain} if $\prec$ is auxiliary, approximating, interpolative and $\vee$-preserving.
\end{dfn}

The `predomain' terminology is due to Keimel \textendash\, see \cite{Keimel2016}.  We immediately see that $\Subset$ is $\cup$-preserving on any $\cup$-basis.  Thus any $\cup$-basis $S\subseteq\mathcal{O}(X)$ of core compact $X$ is a $\Subset$-distributive $\cup$-predomain (see \autoref{CoreCompactDistributivity} and \autoref{CoreCompactInterpolation}).

\begin{prp}
If $S$ is a $\prec$-distributive $\vee$-predomain then, for any $p,q\in S$,
\begin{equation}\label{prec=>Subset'}
p\prec q\qquad\Leftrightarrow\qquad\widehat{S}_p\Subset\widehat{S}_q.
\end{equation}
\end{prp}

\begin{proof}
The $\Rightarrow$ part is immediate from \autoref{prec=>CompactContainment}.  Conversely, say $p\not\prec q$.  As $\prec$ preserves $\vee$, $q^\succ$ is an ideal and hence $(\widehat{S}_r)_{r\prec q}$ generates an ideal of open sets covering $\widehat{S}_q$.  However, none of these open sets covers $\widehat{S}_p$ as $\widehat{S}_p\subseteq\widehat{S}_r$ would imply $p\leq r\prec q$, by \eqref{leq=>subseteq'} and \eqref{prec=>Subset'}.  Thus $\widehat{S}_p\not\Subset\widehat{S}_q$.
\end{proof}

These results can be summarised as giving us a duality of the following classes.
\begin{align*}
\mathbf{{}_\cup B}&=\{S:S\text{ is a $\cup$-basis of a core compact sober space}\}.\\
\mathbf{{}_\prec D_\vee P}&=\{S:S\text{ is a $\prec$-distributive $\vee$-predomain}\}.
\end{align*}

\begin{thm}
$\mathbf{{}_\cup B}$ is dual to $\mathbf{{}_\prec D_\vee P}$.
\end{thm}

More precisely, $\mathbf{{}_\cup B}\subseteq\mathbf{{}_\prec D_\vee P}$, by \autoref{CoreCompactDistributivity} and \autoref{CoreCompactInterpolation}, i.e. any $\cup$-basis of a core compact sober space is a $\Subset$-distributive $\cup$-predomain.  Moreover, the spectrum of such a $\cup$-basis recovers the original space, by \autoref{SoberRecovery}.  Conversely, every $\prec$-distributive $\vee$-predomain $S\in\mathbf{{}_\prec D_\vee P}$ has a core compact sober spectrum $\widehat{S}$, by \autoref{RoundSpectrumSober} and \autoref{CoreCompactSpectrum}.  Moreover, $\leq$ and $\prec$ become $\subseteq$ and $\Subset$ on $(\widehat{S}_p)_{p\in S}\in{}_\cup\mathbf{B}$, by \eqref{leq=>subseteq'} and \eqref{prec=>Subset'}, thus yielding a $\cup$-basis isomorphic to $S$.

\begin{rmk}
The above construction gives us a way of embedding any $\prec$-distributive $\vee$-predomain $S$ into a continuous frame, all we have to do is note that $S$ is isomorphic to $(\widehat{S}_p)_{p\in S}$ which is a subset of the continuous frame of all open sets $\mathcal{O}(\widehat{S})$.  Alternatively, the continuous frame could be constructed directly from the round ideals \textendash\, see \cite[Proposition 5.1.33]{Goubault2013} or \cite[Exercise III.4.17]{GierzHofmannKeimelLawsonMisloveScott2003} \textendash\, it would just have to be verified that $\prec$-distributivity in $S$ yields infinite distributivity of the round ideals.  In \cite{HofmannLawson1978}, Hofmann-Lawson already showed that continuous frames are dual to core compact sober spaces, so this would also provide a somewhat indirect way of obtaining our duality.  Conversely, Hofmann-Lawson duality could be viewed as a special case of our duality when $S$ is a complete lattice on which $\prec$ is defined as the way-below relation.
\end{rmk}

\begin{rmk}
Gr\"atzer's generalisation of Stone duality in \cite[II.5]{Gratzer1978} is also a special case of our duality when $\prec\ =\ \leq$.  Incidentally, as we are already working with bases, one might think it would be harmless to deal with lattices rather than semilattices.  However, taking $\prec\ =\ \leq$ in this case would force the resulting spaces to be coherent and we would only recover Stone's duality for distributive lattices, rather than Gr\"atzer's duality for distributive semilattices.
\end{rmk}

\begin{rmk}We also get Jung-S\"underhauf's duality in \cite{JungSunderhauf1995} as a special case when $\prec$ satisfies their strong proximity axioms.  This points to a potential application of our duality \textendash\, just as proximity relations on subsets of a given space are often used to obtain Tychonoff compactifications, distributive predomain relations could be used to obtain more general sober (even local) compactifications.
\end{rmk}

Another important way of representing bounded distributive lattices is on Priestley spaces (see \cite{Priestley1970}), which are certain kinds of ordered Stone spaces.  It would seem reasonable to guess that $\prec$-distributive $\vee$-predomains could also be represented in a similar way on some class of (potentially connected) pospaces.

\begin{qst}
Does Priestley duality extend to $\prec$-distributive $\vee$-predomains?
\end{qst}

\section{Morphisms}\label{MorphismsS}

Note $\mathbf{{}_\cup B}$ becomes a category when we take partial continuous functions as morphisms, i.e. given $S,S'\in\mathbf{{}_\cup B}$, a morphism from $S$ to $S'$ is a function $\phi$ on (necessarily open) $\mathrm{dom}(\phi)$ such that, for all $p'\in S'$,
\[\phi^{-1}[p']=\bigcup\{p\in S:p\subseteq\phi^{-1}[p']\}.\]
Then $\phi$ yields a relation $\sqsubset_\phi\ \subseteq S\times S'$ given by
\[p\sqsubset_\phi p'\qquad\Leftrightarrow\qquad p\Subset\phi^{-1}[p'].\]
Our goal is to axiomatise such relations in order to make our duality functorial.  At first, we will consider axioms that apply even when we replace $\Subset$ with $\subseteq$, i.e.
\[p\sqsubset p'\qquad\Leftrightarrow\qquad p\subseteq\phi^{-1}[p'].\]

\begin{dfn}\label{Morphisms}
Given $S,S'\in\mathbf{{}_\prec D_\vee P}$, we call $\sqsubset\ \subseteq S\times S'$ a \emph{morphism} if
\begin{align}
\tag{Faithful}\label{Faithful}p\sqsubset0'\qquad&\Rightarrow\qquad p=0.\\
\tag{Auxiliarity}\label{Auxiliarity'}p\leq q\sqsubset q'\leq p'\qquad&\Rightarrow\qquad p\sqsubset p'.\\
\tag{Pushforward}\label{Pushforward}p\prec q\sqsubset r',s'\qquad&\Rightarrow\qquad\exists p'\in S'\ (p\sqsubset q'\prec r',s').\\
\tag{$\vee$-Pullback}\label{veePullback}p\prec q\sqsubset r'\vee s'\qquad&\Rightarrow\qquad\exists q\sqsubset q'\ \exists r\sqsubset r'\ (p\prec r'\vee s').
\end{align}
\end{dfn}

We consider the usual composition of relations, i.e. if $\sqsubset\ \subseteq S\times S'$ and $\sqsubset'\ \subseteq S'\times S''$,
\[p\sqsubset\circ\sqsubset'p''\qquad\Leftrightarrow\qquad\exists p'\in S'\ (p\sqsubset p'\sqsubset p'').\]

\begin{prp}
$\mathbf{{}_\prec D_\vee P}$ forms a category.
\end{prp}

\begin{proof}
It follows from \eqref{Auxiliarity'} that each $S\in\mathbf{{}_\prec D_\vee P}$ has an identity morphism, namely $\leq$.  We just have to show that $\sqsubset\circ\sqsubset'$ is a morphism whenever $\sqsubset$ and $\sqsubset'$ are.  We verify the required properties as follows.
\begin{itemize}
\item[\eqref{Faithful}] If $p\sqsubset p'\sqsubset'0''$ then $p'=0'$ and hence $p=0$, by \eqref{Faithful} for $\sqsubset$ and $\sqsubset'$.
\item[\eqref{Auxiliarity'}] If $p\leq q\sqsubset p'\sqsubset'q''\leq p''$ then $p\sqsubset p'\sqsubset'p''$, by \eqref{Auxiliarity'} for $\sqsubset$ and $\sqsubset'$.
\item[\eqref{Pushforward}] If $p\prec q\sqsubset q'\sqsubset'r'',s''$ then \eqref{Pushforward} for $\sqsubset$ yields $p'\in S'$ with $p\sqsubset p'\prec q'\sqsubset'r'',s''$ and then \eqref{Pushforward} for $\sqsubset'$ yields $q''\in S''$ with $p\sqsubset p'\sqsubset q''\prec r'',s''$.
\item[\eqref{veePullback}] If $p\prec q\sqsubset q'\sqsubset'r''\vee s''$ then \eqref{veePullback} for $\sqsubset'$ yields $r'\sqsubset'r''$ and $s'\sqsubset's''$ with $p\prec q\sqsubset q'\prec r'\vee s'$ and hence $p\prec q\sqsubset r'\vee s'$, by \eqref{Auxiliarity'} for $\sqsubset$.  Then \eqref{veePullback} for $\sqsubset$ yields $r\sqsubset r'$ and $s\sqsubset s'$ with $p\prec r\vee s$. \qedhere
\end{itemize}
\end{proof}

Given $\sqsubset\ \subseteq S\times S'$ define $\phi_\sqsubset$ from $\widehat{S}$ to $\widehat{S'}$ by
\[\phi_\sqsubset(P)=P^\sqsubset\qquad\text{on}\qquad\mathrm{dom}(\phi)=\{P\in\widehat{S}:P^\sqsubset\neq\emptyset\}.\]

\begin{prp}
We have a functor from $\mathbf{{}_\prec D_\vee P}$ to $\mathbf{{}_\cup B}$ given by
\[S\mapsto(\widehat{S}_p)_{p\in S}\qquad\text{and}\qquad\sqsubset\ \mapsto\phi_\sqsubset.\]
\end{prp}

\begin{proof}
For any $S\in\mathbf{{}_\prec D_\vee P}$, we immediately see that $\phi_\leq$ is the identity map on $\widehat{S}$.

The main task is to show $\mathrm{ran}(\phi_\sqsubset)\subseteq\widehat{S'}$, for any morphism $\sqsubset$ between $S,S'\in\mathbf{{}_\prec D_\vee P}$, i.e. that $P^\sqsubset\in\widehat{S'}$ whenever $P\in\widehat{S}$ and $P^\sqsubset\neq\emptyset$.  By \eqref{Faithful}, $P^\sqsubset$ is proper.  Whenever $P\ni p\sqsubset p'\leq q'$, \eqref{Auxiliarity'} yields $p\sqsubset q'$, showing that $P^\sqsubset$ is an up-set.  Whenever $P\ni p\sqsubset r',s'$, the roundness of $P$ yields $t\in P$ with $t\prec p$ and then \eqref{Pushforward} yields $p'\in S'$ with $t\sqsubset p'\prec r',s'$, showing that $P^\sqsubset$ is a $\prec$-filter.  Whenever $P\ni p\sqsubset r'\vee s'$, the roundness of $P$ again yields $t\in P$ with $t\prec p$ and then \eqref{veePullback} yields $r\sqsubset r'$ and $s\sqsubset s'$ with $t\prec r\vee s$.  As $P$ is prime, either $r\in P$ and hence $r'\in P^\sqsubset$ or $s\in P$ and hence $s'\in P^\sqsubset$, showing that $P^\sqsubset$ is prime and hence $P^\sqsubset\in\widehat{S'}$.  Continuity then follows as, for all $p'\in S'$,
\[\phi_\sqsubset^{-1}[\widehat{S'}_{p'}]=\bigcup_{p\sqsubset p'}\widehat{S}_p.\]
Also composition is preserved as, for all $P\in\mathrm{dom}(\phi_{\sqsubset\circ\sqsubset'})=\phi_\sqsubset^{-1}[\mathrm{dom}(\phi_{\sqsubset'})]$,
\[\phi_{\sqsubset\circ\sqsubset'}(P)=P^{\sqsubset\circ\sqsubset'}=(P^\sqsubset)^{\sqsubset'}=\phi_{\sqsubset'}(\phi_\sqsubset(P))=\phi_{\sqsubset'}\circ\phi_\sqsubset(P).\qedhere\]
\end{proof}

On the other hand, the map $\phi\mapsto\ \sqsubset_\phi$ is not (part of) a functor to $\mathbf{{}_\prec D{}_\vee P}$, as it does not preserve identity morphisms.  However, it does preserve composition.

\begin{prp}
For any morphisms $\phi$ from $S$ to $S'$ and $\phi'$ from $S'$ to $S''$ in $\mathbf{{}_\cup B}$,
\[\sqsubset_{\phi'\circ\phi}\ =\ \sqsubset_\phi\circ\sqsubset_{\phi'}.\]
\end{prp}

\begin{proof}
Say $p\sqsubset_{\phi'\circ\phi}p''$, which means $p\Subset\phi^{-1}[\phi'^{-1}[p'']]$ and hence we have compact $C$ with $p\subseteq C\subseteq\phi^{-1}[\phi'^{-1}[p'']]$.  Thus $\phi[C]$ is also compact and $\phi[C]\subseteq\phi'^{-1}[p'']$.  As $S'$ is a $\cup$-basis, this means we have $p'\in S'$ with $\phi[C]\subseteq p'\Subset\phi'^{-1}[p'']$ and hence $p\sqsubset_\phi p'\sqsubset_{\phi'}p''$.  The converse is immediate from the definitions.
\end{proof}

\begin{prp}
For any morphism $\phi$ from $S$ to $S'$ in $\mathbf{{}_\cup B}$, we always have
\[S'_{\phi(x)}=\phi_{\sqsubset_\phi}(S_x).\]
\end{prp}

\begin{proof}
Just note that $p'\in S'_{\phi(x)}$ iff, for some $p\in S$, $x\in p\Subset\phi^{-1}[p']$ which means $S_x\ni p\sqsubset_\phi p'$, i.e. $p'\in\phi_{\sqsubset_\phi}(S_x)$.
\end{proof}

For any $\sqsubset\ \subseteq S\times S'$, we define $\sqsubset_\vee\ \subseteq S\times S'$ by
\[p\sqsubset_\vee p'\qquad\Leftrightarrow\qquad\exists\text{ finite }F\sqsubset p'\ (p\prec\bigvee F).\]

\begin{prp}
For any morphism $\sqsubset$ between $S,S'\in\mathbf{{}_\prec D{}_\vee P}$,
\[p\sqsubset_\vee p'\qquad\Leftrightarrow\qquad\widehat{S}_p\sqsubset_{\phi_\sqsubset}\widehat{S'}_{p'}.\]
\end{prp}

\begin{proof}
Say $p\sqsubset_\vee p'$, so we have finite $F\sqsubset p'$ with $p\prec\bigvee F$ and hence $\widehat{S}_p\Subset\widehat{S}_{\bigvee F}$, by \eqref{prec=>Subset'}.  But $F\sqsubset p'$ implies $\widehat{S}_{\bigvee{F}}=\bigcup_{f\in F}\widehat{S}_f\subseteq\phi_\sqsubset^{-1}[\widehat{S'}_{p'}]$ and hence $\widehat{S}_p\Subset\phi_\sqsubset^{-1}[\widehat{S'}_{p'}]$, i.e. $\widehat{S}_p\sqsubset_{\phi_\sqsubset}\widehat{S'}_{p'}$.

Conversely, say $\widehat{S}_p\Subset\phi_\sqsubset^{-1}[\widehat{S'}_{p'}]$.  Take $q\in S$ with $\widehat{S}_p\Subset\widehat{S}_q\Subset\phi_\sqsubset^{-1}[\widehat{S'}_{p'}]$.   The definition of $\phi_\sqsubset$ means we can cover $\phi_\sqsubset^{-1}[\widehat{S'}_{p'}]$ by sets $\widehat{S}_f$ such that $f\sqsubset p'$.  Then the definition of $\Subset$ means we have finite $F\sqsubset p'$ with $\widehat{S}_p\Subset\widehat{S}_q\subseteq\bigcup_{f\in F}\widehat{S}_f=\widehat{S}_{\bigvee F}$ and hence $p\prec\bigvee F$, by \eqref{prec=>Subset'}, i.e. $p\sqsubset_\vee p'$.
\end{proof}

These results show that we can obtain a category equivalent to $\mathbf{{}_\cup B}$ if we restrict to `$\vee$-morphisms' in $\mathbf{{}_\prec D{}_\vee P}$, i.e. those morphisms with $\sqsubset\ =\ \sqsubset_\vee$.  Equivalently, these are the morphisms satisfying the following extra conditions.
\begin{align}
\tag{Left Interpolation}p\sqsubset p'\qquad&\Rightarrow\qquad\exists q\in S\ (p\prec q\sqsubset p').\\
\tag{$\vee$-Preserving}q,r\sqsubset p'\qquad&\Rightarrow\qquad q\vee r\sqsubset p'.
\end{align}
Note each $S$ still has an identity in this new category, it is just $\prec$ instead of $\leq$.  Incidentally, if we wanted to restrict to total functions in $\mathbf{{}_\cup B}$, then we would add one further condition on the morphisms in $\mathbf{{}_\prec D{}_\vee P}$, namely
\[\hspace{77pt}\tag{Total}p\prec q\qquad\Rightarrow\qquad\exists q'\in S'\ (p\sqsubset q').\]


\bibliography{maths}{}
\bibliographystyle{alphaurl}

\end{document}